\documentclass[11pt]{amsart}
\usepackage{amssymb}
\usepackage{amsmath}
\usepackage[all]{xy}
\usepackage{amsthm}
\usepackage{amscd} 
\newtheorem{thm}{Theorem}[section]
\newtheorem{rk}[thm]{Remark}
\newtheorem{prop}[thm]{Proposition}   
\newtheorem{clly}[thm]{Corollary}
\newtheorem{lemma}[thm]{Lemma}
\newtheorem{defi}[thm]{Definition}

\newtheorem{nota}[thm]{Notation}

\newcommand{\bpf}{\begin{proof}}
\newcommand{\epf}{\end{proof}}
\newcommand{\bprop}{\begin{prop}}
\newcommand{\eprop}{\end{prop}}
\newcommand{\bthm}{\begin{thm}}
\newcommand{\ethm}{\end{thm}}
\newcommand{\brk}{\begin{rk}}
\newcommand{\erk}{\end{rk}}
\newcommand{\bdefi}{\begin{defi}}
\newcommand{\edefi}{\end{defi}}
\newcommand{\blemma}{\begin{lemma}}
\newcommand{\elemma}{\end{lemma}}
\newcommand{\bclly}{\begin{clly}}
\newcommand{\eclly}{\end{clly}}
\newcommand{\bnota}{\begin{nota}}
\newcommand{\enota}{\end{nota}}
\newcommand{\sss}{\scriptscriptstyle}

\newcommand{\kex}{{\mathcal {K}}(E_X)}

\newcommand{\id}{\mbox{id}}
\newcommand{\im}{\mbox{Im}}

\newcommand{\tx}{\tilde{x}}

\newcommand{\aut}{\mbox{Aut}}

\newcommand{\cstar}{\mbox{$C^*$}}
\newcommand{\acx}{A\rtimes X}
\newcommand{\acxcg}{(A\rtimes X)\rtimes_{\alpha}G}

\newcommand{\CC}{\mathbb{C}}

\newcommand{\ZZ}{\mathbb{Z}}
\newcommand{\supp}{\mbox{supp}}
\newcommand{\cl}{\mathcal{L}}
\newcommand{\ck}{\mathcal{K}}

\newcommand{\lran}{\rangle_{{\scriptscriptstyle L}}}
\newcommand{\rran}{\rangle_{{\scriptscriptstyle R}}}
\newcommand{\tnx}{T^{n}_x}
\newcommand{\xn}{X^{\otimes n}}
\newcommand{\tX}{{\tilde X}}

\title{Takai duality for crossed products by Hilbert \cstar-bimodules}
\author{Beatriz Abadie} 
\thanks{Partially supported by Proyecto Fondo Clemente Estable 8013}
\address{Centro de Matem\'atica. Facultad de Ciencias. Igu\'a 4225, CP 11 400, Montevideo, Uruguay.}
\email{abadie@cmat.edu.uy}

\begin{document}

\begin{abstract}
 We discuss the crossed product by the dual action $\delta$ of the circle  on the crossed product $\acx$  of a \cstar-algebra $A$ by a Hilbert \cstar-bimodule $X$. When $X$ is an $A-A$  Morita equivalence bimodule,  the double crossed product $A\rtimes X\rtimes_\delta S^1$ is shown to be Morita equivalent to the \cstar-algebra $A$.
\end{abstract}

\maketitle
\section{Introduction}

The crossed product $A\rtimes X$ of a \cstar-algebra $A$ by a Hilbert \cstar-bimodule $X$ was introduced in \cite{aee} and shown to be a generalization of the crossed product by an automorphism. There is an obvious generalization of the dual action to this context, which raises the question of whether there is an analog of Takai duality (\cite{takai}). 
We show in this work that when $X$ is an $A-A$ Morita equivalence bimodule, that is, when it is a full Hilbert \cstar-module both on the left and the right, then the double crossed product $A\rtimes X\rtimes_\delta \ZZ$ is Morita equivalent to the \cstar-algebra $A$. Namely, if $E_X$ denotes the right  Hilbert \cstar-module over $A$ defined by $E_X=\oplus_{n\in \ZZ}X^{\otimes n}$, then we identify $A\rtimes X\rtimes_\delta\ZZ$  with $\mathcal {K}(E_X)$, the \cstar-algebra of compact operators on $E_X$, and we describe the double dual action on $\mathcal {K}(E_X)$.  Our proof heavily relies on the  universal properties of the crossed products by  an automorphism and by a Hilbert \cstar-bimodule, much as in \cite{raeb}.

This work is organized as follows. After establishing some preliminary results and notation  in section 2, we introduce in section 3  representations on the crossed product $A\rtimes X$ induced by representations on $A$. Section 4 is devoted to  the discussion of certain actions of amenable locally compact groups on $\acx$ that leave $A$ and $X$ invariant. We show that the crossed product of $\acx$ by an action of this kind can be written as the crossed product of  a \cstar-algebra by a Hilbert \cstar-bimodule. 
These results enable us to  represent, in section 5, the double crossed product $A\rtimes X\rtimes_\delta\ZZ$ as adjointable operators on $E_X$. When $X$ is a Morita equivalence bimodule this representation turns out to be an isomorphism onto $\mathcal {K}(E_X)$. This yields the Morita equivalence between $A\rtimes X\rtimes_\delta\ZZ$ and $A$.

\vspace{.1in}
\noindent {\bf{Acknowledgement.}} Part of this research was carried out during my visit to the  Universit\'{e} d'Orl\'eans. I would like to   express my gratitude to  Jean Renault for his kind invitation and to   the members  of the  department of mathematics for their warm hospitality.

\section{Preliminaries}
\label{prel}

We next establish our basic notation concerning Hilbert \cstar-modules and bimodules. We refer the reader to \cite{la} for further details.

Let $X$ and $Y$ be  right Hilbert \cstar-modules over a \cstar-algebra $A$.  We denote by  $ \cl(X, Y)$   the space of adjointable maps from $X$ to $Y$ and by $ {\mathcal K}(X,Y)$ the space of compact operators, that is, the closed subspace  spanned by $\{ \theta_{y,x}:\  x\in X,\ y\in Y\}$, where $\theta_{y,x}:X\rightarrow Y$ is given by $\theta_{y,x}(z)=y\langle x,z \rangle$. We  will also use  the notation above when $X$ and $Y$ are Hilbert \cstar-bimodules, thus viewing them as right Hilbert modules. Undecorated inner products will always denote right inner products.

Throughout this  work we consider Hilbert \cstar-bimodules in the sense of \cite{bms}. That is,  a Hilbert \cstar-bimodule $X$ over a \cstar -algebra $A$ consists of a vector space $X$ which is both a right and a left Hilbert \cstar-module over $A$ and satisfies  $\langle x,y\lran z =x\langle y,z\rran$ and $ (ax)b=a(xb)$, for all $x,y,z\in X$  and $a,b\in A$.

Let $X$ and $Y$ be Hilbert \cstar-bimodules over the \cstar-algebras $A$ and $B$, respectively. A morphism of Hilbert \cstar-bimodules 
\[(\phi_A,\phi_X):(A,X)\rightarrow (B,Y)\] consists of a $*$-homomorphism  $\phi_A:A\rightarrow B$ and a linear map $\phi_X:X\rightarrow Y$ such that
\begin{gather*}
\phi_X(ax)=\phi_A(a)\phi_X(x),\  \phi_X(xa)=\phi_X(x)\phi_A(a), \\
\langle \phi_X(x),\phi_X(y)\rangle_{{\scriptscriptstyle L}}=\phi_A(\langle x,y\rangle_{{\scriptscriptstyle L}}),\mbox{  and  }\langle \phi_X(x),\phi_X(y)\rangle_{{\scriptscriptstyle R}}=\phi_A(\langle x,y\rangle_{{\scriptscriptstyle R}}),
\end{gather*}
for all $x,y\in X$ and $a\in A$.
\brk
\label{morph}
{\rm The map $\phi_X$  above is norm decreasing, and isometric when so is $\phi_A$, since
\[\|\phi_X(x)\|^2=\|\langle \phi_X(x),\phi_X(x)\rangle\|=\|\phi_A(\langle x,x \rangle)\|\leq \|\langle x,x \rangle \|=\|x\|^2.\]
}
\erk
A  representation of $(A,X)$ on a \cstar-algebra $B$ consists of a  morphism $(\phi_A,\phi_X):(A,X)\rightarrow (B,B)$  where   $B$ is viewed as a Hilbert \cstar-bimodule over itself in   the usual way.

The crossed product $\acx$  is (\cite{aee}) the universal \cstar-algebra carrying  a  representation $(i_A,i_X)$ of $(A,X)$. Besides, the maps $i_A$ and $i_X$  are isometric (\cite[2.10]{aee}), and $\acx$ is generated as a \cstar-algebra by the images of $i_A$ and $i_X$.
The crossed product $\acx$ carries an action of $S^1$, called the dual action, which is the identity on $i_A(A)$ and is given by $i_X(x)\mapsto \lambda i_X(x)$ for $x\in X$ and $\lambda\in S^1$. Moreover, $i_A(A)$ and $i_X(X)$ are the fixed-point subalgebra and the first spectral subspace of the dual action, respectively  (see \cite[3]{aee}).

\begin{rk}
\label{covmaps} 
Let $(\phi_A,\phi_X):(A,X)\longrightarrow (B,Y)$ be a morphism of Hilbert \cstar-bimodules. Then there is a unique $*$-homomorphism 
\[\phi_A\rtimes \phi_X: A\rtimes X\longrightarrow B\rtimes Y\] such that the diagram
\[\begin{CD}
A\rtimes X @>\phi_A\rtimes \phi_X>> B\rtimes Y\\
@A(i_{A}, i_X)AA             @AA(i_{B}, i_Y)A\\
(A,X) @>(\phi_A, \phi_X)>>(B,Y)
\end{CD}\]     
commutes. The map $\phi_A\rtimes \phi_X$ is injective when so is $\phi_A$, and it is surjective if so are $\phi_A$ and $\phi_X$.

Besides, the correspondence  $(A,X)\mapsto \acx,\ (\phi_A,\phi_X)\mapsto \phi_A\rtimes \phi_X$ is functorial.
\end{rk}

\begin{proof}
Since  $(i_B\circ \phi_A, i_Y\circ \phi_X)$ is a  representation of $(A,X)$ on $B\rtimes Y$,  the existence and the uniqueness of $\phi_A\rtimes\phi_X: A\rtimes X\longrightarrow B\rtimes Y$ follow from the universal property of $\acx$. 

The map $\phi_A\rtimes\phi_X$ is covariant for the dual actions of $S^1$ on $A\rtimes X$ and $B\rtimes Y$ respectively, so it is injective if and only if it is injective when restricted to the fixed-point algebra (\cite[2.9]{cact}),  that is, when $\phi_A$ is injective. When $\phi_A$ and $\phi_X$ are surjective, the image of $\phi_A\rtimes\phi_X$ contains $B$ and $Y$, so it is all of $B\rtimes Y$. The last statement is apparent.

 \end{proof}

When   $X$ and $Y$ are Hilbert \cstar-bimodules over $A$ and $x\in X$, we denote by    $T_x\in  \cl(Y,X\otimes Y)$ the creation operator    defined by $T_x(y)=x\otimes y$, for $y\in Y$. 

The following facts are well known and easy to check: 
\begin{gather}
T_{ax}(y)=aT_x(y),\  T_{xa}(y)=T_x(ay), \label{E:creation1}\\ 
\|T_x\|\leq \|x\|,\  T^*_{x_0} (x\otimes y)=\langle x_0,x\rran y,\label{E:creation2}\\
 T_{x_0}T^*_{x_1}(x\otimes y)=\langle x_0,x_1\lran x\otimes y  \mbox{ and } T^*_{x_0}T_{x_1}(y)=\langle x_0,x_1\rran y \label{E:creation3},
\end{gather}
for $a\in A$,  $x,x_0,x_1\in X$ and $y\in Y$.

\begin{defi}
\label{TnEx} {\rm Let $X$ be a Hilbert \cstar-bimodule over the \cstar-algebra $A$,  $x\in X$, and let $n$ be a non-negative integer.  We denote by $ T^{n}_x$ the map $ T^{n}_x \in \cl\left(\xn,X^{\otimes n+1}\right)$ described above, where   $X$ and $X\otimes A$ are identified in the usual way when $n=0$.

When  $n<0$, $X^{\otimes n}$ denotes $(\tX)^{\otimes -n}$,  $\tX$ being the dual bimodule of $X$ as defined in \cite{irep}. That is, $\tilde{X}$ is the conjugate vector space of $X$   and carries the $A$-Hilbert \cstar-bimodule structure given by
\begin{gather*}
a\cdot \tilde{x}=\widetilde{xa^*}, \tilde{x}\cdot a= \widetilde{a^*x},\ 
\langle \tilde{x}, \tilde{y}\lran = \langle x,y\rran,\ \langle \tilde{x}, \tilde{y} \rran = \langle x,y\lran,
\end{gather*}
where $\tilde x$ denotes the element $x\in X$ viewed as an element of the dual bimodule ${\tilde X}$. For negative values of $n$ we define
 \begin{align}
 \tnx:= \big(T^{-n-1}_{{\tilde x}}\big)^*.\label{E:tnegn} 
\end{align}
Thus, for $n<0$, 
\begin{gather}\label{E:negn} \tnx(\tx_1\otimes\tx_2\otimes\dots\otimes \tx_{-n})=\langle x,x_1\lran\ \tx_2\otimes\dots\otimes\tx_{-n}\\
\label{E:negnstar} (T_x^n)^*(\tx_1\otimes\tx_2\otimes\dots\otimes \tx_{-n-1})= \tx\otimes\tx_1\otimes\tx_2\otimes\dots\otimes\tx_{-n-1}.
\end{gather}
Let $L^n$   denote the left action of $A$ on $X^{\otimes n}$ for $n\in \ZZ$. It follows from equations (\ref{E:creation1}) - (\ref{E:creation3}) for $n\geq 0$,  and from some straightforward maneuvering and equations (\ref{E:negn}) and  (\ref{E:negnstar}) for $n<0$,  that
\begin{gather}
T^n_{ax}=L^{n+1}_aT^n_x,\ T^n_{xa}=T^n_xL^{n}_a. \label{E:tn1}\\
(T^n_x)(T^n_y)^*=L^{n+1}_{\langle x,y\lran},\ (T^n_x)^*(T^n_y)=L^n_{\langle x,y\rran},  \label{E:tn2}
\end{gather}
for all $n\in \ZZ$, $a\in A$, and $x,y\in X$.

Throughout this work we will denote by  $E_X$  the right $A$-Hilbert \cstar-module defined by $E_X=\bigoplus_{-\infty}^{+\infty} \xn$. 

\begin{rk}
\label{included}
{\rm  We will often view $\cl(X^{\otimes n}, X^{\otimes m})$ as a closed subspace of $\cl(E_X)$ by means of the isometric linear map $i_{n,m}:\cl(X^{\otimes n}, X^{\otimes m})\rightarrow \cl(E_X)$ given by 
\[ [(i_{n,m}T)(\eta)](k)=\delta_m(k)T(\eta(n)), \]
for $T\in \cl(X^{\otimes n}, X^{\otimes m})$, $\eta\in E_X$, and $k\in \ZZ$.

Note that 
\[i_{n,m}(T)^* =i_{m,n}(T^*),\mbox{ and }  i_{m,p}(S)\circ i_{n,m}(T)=i_{n,p}(S\circ T).\]
Besides, 
\[ i_{n,m}(\theta_{u,v})=  \theta_{u\delta_m,v\delta_n},\]
so under this identification $\ck(X^{\otimes n}, X^{\otimes m})\subset \ck(E_X)$.}
\end{rk}

\section{Induced representations on crossed products}

We discuss in this section certain representations of $\acx$ that are induced from representations of $A$.

Let $X$  be a Hilbert \cstar-bimodule over $A$, and let $E_X$ be  the right Hilbert \cstar-module defined at the end of section 2. 

We define  $(\Lambda_A, \Lambda_X):(A,X)\rightarrow \cl(E_X)$ by 
\begin{gather}
[\Lambda_A(a)(\eta)](n)=L^n_a(\eta(n)), \quad [\Lambda_X(x)(\eta)](n)=T^{(n-1)}_x(\eta(n-1)),\label{E:lambdas}
\end{gather}
for $a\in A$, $x\in X$,  $\eta\in E_X$, and $T^{(n)}_x$ as in Definition \ref{TnEx}.

It is easily checked that $(\Lambda_X(x)^*(\eta))(n)=(T^n_x)^*(\eta(n+1))$. It follows from equations (\ref{E:tn1}) and (\ref{E:tn2}) that $(\Lambda_A, \Lambda_X)$ is a representation on $\cl(E_X)$ and  induces a *-homomorphism 
\[\Lambda:=\Lambda_A\times\Lambda_X:\acx\longrightarrow \cl(E_X),\]
as in Remark \ref{covmaps}
}
\end{defi}
\begin{prop}
\label{action} Let $A$, $X$ and $\Lambda :\acx\longrightarrow \cl(E_X) $ be as  above. For $\lambda\in S^1$, let $U_{\lambda}$ be the unitary operator on $E_X$ defined by 
\[(U_{\lambda}\eta)(n)=\lambda^n\eta(n),\]
 for $\eta\in E_X,\ n\in \ZZ$. Then:
\begin{enumerate}

\item Conjugation by $U_\lambda$ defines a strongly continuous action $\beta$ of $S^1$ on the image of $\Lambda$.
\item The map  $\Lambda$  is injective and covariant for the dual action and $\beta$.
\end{enumerate}
\end{prop}
\begin{proof}
It is   easily checked that the map  $\lambda\mapsto U_{\lambda}$ is a group homomorphism from $S^1$ to the unitary group of $\cl(E_X)$  and that
\begin{gather}\label{E:ulambda}
 U_{\lambda}\Lambda_A(a)U^*_{\lambda}=\Lambda_A(a),\  U_{\lambda}\Lambda_X(x)U^*_{\lambda}=\lambda \Lambda_X(x), 
\end{gather}
for all $x\in X$  and $a\in A$. 

On the other hand, the set
\[\{T\in \cl(E_X):\ \lambda\mapsto U_\lambda TU^*_\lambda \mbox{ is continuous on } S^1\}\]
is a \cstar-algebra of $\cl(E_X)$. Therefore, since $\acx$ is generated as a \cstar-algebra by $i_A(A)$ and $i_X(X)$, conjugation by $U_{\lambda}$ defines a strongly continuous action $\beta$ of $S^1$ on the image of $\Lambda$. An analogous reasoning, together with equation (\ref{E:ulambda}), shows that $\Lambda$ is covariant for $\beta$ and the dual action on $\acx$. Finally, since  the restriction $\Lambda_A$ of $\Lambda$ to the fixed-point subalgebra is injective,  so is $\Lambda$ by \cite[2.9]{cact}. 
\end{proof}
\begin{rk} {\rm If $\acx$ is viewed as the cross-sectional \cstar-algebra of a Fell bundle  as in \cite[2.9]{aee}, then $\Lambda$  is, in the terminology of \cite{amen}, the left regular representation and its injectivity follows from \cite{amen} and the amenability of $\ZZ$.}
\end{rk}

\begin{defi}
\label{indrep}
{\rm{Let $X$ be a Hilbert \cstar-bimodule over a  \cstar-algebra A. A non-degenerate representation $\pi$ of $A$ on a Hilbert space $H$ gives rise to the representation $\Lambda\otimes\id_H$ of $\acx$ on the Hilbert space $E_X\otimes_{\pi}H$. We will refer to $\Lambda\otimes\id_H$ as the  representation of $\acx$ induced by $\pi$.}}
\end{defi}
\brk
\label{aalpha}
{\rm{ For $\alpha\in \aut (A)$ let $A_\alpha$ be the $A$-Hilbert \cstar-bimodule
consisting of $A$ as a vector space with structure defined by
\[a\cdot x=ax,\ x\cdot a=x\alpha(a),\ \langle x,y\lran =xy^*\mbox{ and } \langle x,y\rran=\alpha^{-1}(x^*y),\]
for $a\in A$ and $x,y\in A_\alpha$.
It was shown in \cite[3.2]{aee} that there is an isomorphism  $J:A\rtimes A_\alpha\rightarrow A\rtimes_\alpha \ZZ$  given by 
\begin{gather}
J(i_A(a)) =a\delta_0\in C_c(\ZZ, A)\mbox{  and }J(i_X(x)) =x\delta_1\in C_c(\ZZ, A),\label{E:defJ} 
\end{gather}
for $a\in A$ and $x\in A_\alpha$.
Denote by  $I_n:(A_\alpha)^{\otimes n}\rightarrow A,\ n\in \ZZ$
 the map defined by 
\[\left\{
\begin{array}{l}
I_n(a_1\otimes\dots\otimes a_n)=\alpha^{-n}(a_1)\alpha^{-n+1}(a_2)\dots\alpha^{-1}(a_n)\mbox{ for }n\geq 0 \\
I_n=\id_A, \mbox{ for }n= 0 \\
I_n(\tilde{a}_1\otimes\dots\otimes \tilde{a}_{-n})=\alpha^{-n-1}(a^*_1)\alpha^{-n-2}(a^*_2)\dots a^*_{-n},\mbox{ for }n\geq 0 
\end{array}
\right.\]
Straightforward computations show that $I_n$ is a homomorphism of right Hilbert \cstar-modules over $A$ and that
\begin{gather}
I_n(L^n(a)(c))=\alpha^{-n}(a)I_n(c) \label{E:in1}\\
I_{n+1}(T^n_x(c))=\alpha^{-(n+1)}(x)I_{n}(c),\label{E:in2}
\end{gather}
for all $a\in A$, $x\in A_\alpha$ and $c\in (A_{\alpha})^{\otimes n}$.

Now, given a non-degenerate representation $\pi$ of $A$  on a Hilbert space H, we define $U:E_X\otimes_\pi H\rightarrow l^2(\ZZ,H)$  by 
\[ [U(\eta\otimes h)](n)=\pi[I_n(\eta(n))](h),\]
for $ X=A_\alpha,\ \eta\in E_X,\ h\in H$, and $n\in \ZZ$.
Note that  $U$ extends to a unitary operator because 
\[\begin{array}{c}
\langle U(\sum_{i=1}^p \eta_i\otimes h_i), U(\sum_{j=1}^q \xi_j\otimes k_j)\rangle\\
\\
=\sum_{i,j,n}\langle \pi[I_n(\eta_i(n))](h_i), \pi[I_n(\xi_j(n))](k_j)\rangle\\
\\
= \sum_{i,j}\langle h_i,\pi(\langle \eta_i,\xi_j\rran) (k_j)\rangle\\
\\
=\langle \sum_i \eta_i\otimes h_i, \sum_j \xi_j\otimes k_j\rangle ,
\end{array}\]
for $\eta_i,\xi_j\in E_X$ and $ h_i, k_j\in H$.

Let now $\pi_{\alpha}\times \lambda$ denote the representation
of $A\rtimes_\alpha \ZZ$ on $l^2(\ZZ, H)$ induced by $\pi$. That is, $\pi_{\alpha}\times \lambda$ is the integrated form of the covariant pair $(\pi_{\alpha},\lambda)$  defined by 
\[[(\pi_{\alpha}(a))(\xi)](n)= [\pi(\alpha^{-n}(a))](\xi(n)), \ (\lambda_k\xi)(n)=\xi(n-k),\]
for $\xi\in l^2(\ZZ, H)$.
Then  the diagram
\[
\begin{CD}
 E_X\otimes_{\pi}H @>U>> l^2(\ZZ,H) \\
 @V{(\Lambda\otimes \id_H)}(\phi)VV      @VV{(\pi_{\alpha}\times \lambda)J(\phi)}V\\
 E_X\otimes_{\pi}H  @>U>> l^2(\ZZ,H)
\end{CD}
\]
commutes for all $\phi\in A\rtimes A_\alpha$ and $J$ as in equation (\ref{E:defJ}). 

It suffices to check this statement for $\phi=i_A(a)$ and $\phi=i_X(x)$, for $a\in A$ and $x\in A_\alpha$, which follows from equations (\ref{E:in1}) and (\ref{E:in2}) above.}}
\erk
\brk

\label{rkindrep}
{\rm{Let  $X$ be a Hilbert \cstar-bimodule over $A$,  and let $\pi$ be a non-degenerate representation of $A$ on a Hilbert space $H$. Denote by $V:A\otimes_\pi H\rightarrow H$ the unitary operator defined by $V(a\otimes h)=\pi(a)(h)$.
Let $i$ be the isometric embeddding of $A\otimes_\pi H$ into $E_X\otimes_\pi H$ given by $i(a\otimes h)=a\delta_0\otimes h$, and let $S: H\longrightarrow E_X\otimes_{\pi}H$ be the isometry defined  by $S=i\circ V^*$. 

Then  the diagram

\[
\begin{CD}
H  @>S>> E_X\otimes_{\pi}H \\
 @V\pi(a)VV      @VV{(\Lambda\otimes\id_H)(i_A(a))}V\\
 H  @>S>>E_X\otimes_{\pi}H 
\end{CD}
\]
commutes for all $a\in A$.}}
\erk
\begin{proof}
Straightforward computations prove the statement.
\end{proof}
\begin{prop}
Let $X$ be a Hilbert \cstar-bimodule over $A$, and let $\pi$ be a faithful non-degenerate representation of $A$ on a Hilbert space $H$. Then the induced representation $\Lambda\otimes \id_H$ on $ E_X\otimes_{\pi}H$ is faithful.
\end{prop}
\begin{proof}
Let $\beta$ denote the strongly continuous action of $S^1$ on the image of $\Lambda$ defined  in Proposition \ref{action}. Then $\beta\otimes \id$ is a strongly continuous action of $S^1$ on the image of $\Lambda\otimes \id_H$. Besides, by Proposition \ref{action} (2), $\Lambda\otimes \id_H$ is covariant for the dual action $\delta$  and $\beta\otimes\id$.
  
Thus, by \cite[2.9]{cact}, it suffices to show that  $\Lambda\otimes \id$ is injective on the fixed-point subalgebra $i_A(A)$. This last fact follows from Remark \ref{rkindrep} and the injectivity of $\pi$.
\end{proof}

\section{Covariant actions on crossed products}

Throughout this section all integrals over groups are taken  with respect to Haar measure.

\begin{defi} \label{dcovact}Let $G$ be a locally compact group, and let $X$ be a Hilbert \cstar-bimodule over a \cstar-algebra $A$. 
A strongly continuous covariant action $(\alpha_{{\scriptscriptstyle A}},\alpha_{{\scriptscriptstyle X}})$ of $G$ on $(A,X)$ consists of a strongly continuous action $\alpha_{{\scriptscriptstyle A}}$ of $G$ on $A$ and a group homomorphism $\alpha_{{\scriptscriptstyle X}}$ from $ G$ to the group of invertible linear maps on $X$  such that  
\begin{enumerate}
\item The map $t\longrightarrow (\alpha_{{\scriptscriptstyle X}})_t(x)$ is continuous on $G$ for all $x\in X$.
\item The pair  $((\alpha_{{\scriptscriptstyle A}})_t,(\alpha_{{\scriptscriptstyle X}})_t)$ is an isomorphism of Hilbert \cstar-bimodules for all $t\in G$. 
\end{enumerate} 
\begin{rk}
{\rm For $(\alpha_A,\alpha_X)$ as above, the maps $(\alpha_{{\scriptscriptstyle X}})_t$ in Definition \ref{dcovact} are isometric for all $t\in G$, by Remark \ref{morph}. }

\end{rk}
\end{defi}

\brk 
\label{covact}Let $X$ be a Hilbert \cstar-bimodule over a \cstar-algebra $A$ and let $(\alpha_{{\scriptscriptstyle A}},\alpha_{{\scriptscriptstyle X}})$ be a strongly continuous action of a locally compact group $G$ on $(A,X)$.
 Then $\alpha_t:=(\alpha_{{\scriptscriptstyle A}})_t\rtimes(\alpha_{{\scriptscriptstyle X}})_t$ defines a strongly continuous action $\alpha$ of $G$ on $A\rtimes X$.
\erk
\begin{proof}
It follows from Remark \ref{covmaps} that $t\mapsto\alpha_t $ is a group homomorphism from $G$ to $\aut(A\rtimes X)$.
Since $\|\alpha_t\|=1$ for all $t\in G$, the set $\{c\in \acx:\ t\mapsto \alpha_t(c) \mbox{ is continuous}\}$ is a \cstar-subalgebra of $\acx$ containing $i_A(A)$ and $i_X(X)$, which shows that $\alpha$ is strongly continuous.
\end{proof}

\begin{rk}
{\rm  Let  $f$ be in $C_c(G,X)$ for a locally compact group $G$ and  a Hilbert \cstar-bimodule $X$ over a \cstar-algebra $A$. By identifying $X$ with its isometric copy $i_X(X)$ in $\acx$, we view  $f$  as an $\acx$-valued map. Thus the integral $\int_G f d\mu$  has its usual meaning. Notice that this integral belongs to $i_X(X)$, since it can be approximated by sums $\sum_1^nc_if(t_i)$, for $c_i\in \CC$ and $t_i\in G$. This is the way we will view $\int_Gf d\mu$ as an element of $X$ throughout this work.

The same procedure could be followed in order to view   $A$-valued functions  as being  $\acx$-valued. The integral does not depend on the approach, since the restriction to $A$ of a non-degenerate faithful representation of $\acx$ is again a  non-degenerate faithful representation. }
\end{rk}
\begin{prop} 
\label{dcp}
Let $X$ be a Hilbert \cstar-bimodule over $A$, and let $\alpha$ be the strongly continuous action on $A\rtimes X$ induced by a covariant action $(\alpha_{{\scriptscriptstyle A}},\alpha_{{\scriptscriptstyle X}})$ of an amenable  locally compact group $G$ on $(A,X)$. 

Let $i_A$ and $i_X$ be the embeddings of $A$ and $X$, respectively, in $\acx$.
Then the  map that sends $\phi\in C_c(G,i_A(A))$ to $i^{-1}_A\circ \phi$ and $f\in C_c(G,i_X(X))$ to $i^{-1}_X\circ f$ extends to an isomorphism from  $(A\rtimes X)\rtimes_{\alpha}G$ to $ (A\rtimes_{\alpha_{{\scriptscriptstyle A}}}G)\rtimes Y,$ 
where $Y$ is the Hilbert \cstar-bimodule over $A\rtimes_{\alpha_{{\scriptscriptstyle A}}}G$ consisting of  the completion of $C_c(G,X)$ with the structure defined by 
\[(\phi f)(t)=\int_G\phi(u)(\alpha_{\scriptscriptstyle X})_u[f(u^{-1}t)]du,\  (f\phi)(t)=\int_G f(u)(\alpha_{\scriptscriptstyle A})_u[\phi(u^{-1}t)]du,\]
\[\langle f,g\lran(t) =\int_G \Delta(t^{-1}u)\langle f(u),(\alpha_{\scriptscriptstyle X})_t[g(t^{-1}u)]\lran du\ \mbox{and}\ \]
\[ \langle f,g\rran (t)= \int_G (\alpha_{\scriptscriptstyle A })_{u^{-1}}\langle f(u),g(ut)\rran du,\]
where   $f,g\in C_c(G,X), \  \phi\in C_c(G,A)$, and $\Delta$ denotes the modular function on $G$.
 
\end{prop}
 
\begin{proof}

Let $\sigma$ denote the dual action of $S^1$ on $A\rtimes X$.  Note that $\sigma$ and  $\alpha$ commute, since so do their restrictions to the images of $i_A$ and $i_X$.  Therefore (\cite[1.2]{fpa}) $\sigma$ induces an action $\gamma$ of $S^1$ on $(A\rtimes X)\rtimes_{\alpha}G$ given by 
\[[\gamma_{\lambda}(\phi)](t)=\sigma_{\lambda}[\phi(t)],\] 
for $\phi\in C_c(G,A\rtimes X)$, $t\in G$, and $\lambda\in S^1$.

Let $B_0$ and $B_1$ be the closures in $\acxcg$ of the sets of functions in $C_c(G,\acx)$ whose image lies, respectively, in $i_A(A)$ and $i_X(X)$.
We next show that $B_0$ and $B_1$ are  the fixed-point subalgebra and the first spectral subspace, respectively, for the action $\gamma$. We will then prove that the action $\gamma$ is semi-saturated and that $(B_0, B_1)$ and $( A\rtimes_{\alpha_{{\scriptscriptstyle A}}}G, Y)$ are isomorphic as Hilbert \cstar-bimodules. At this point the statement will follow from \cite[3.1 ]{aee}. 

Let  $P_n$ denote   the $n^{\rm{th}}$ spectral projection for the action $\gamma$, that is, 
\[P_n(c)=\int_{S^1}\lambda^{-n}\gamma_{\lambda}(c)d\lambda.\] 
The  $n^{\rm{th}}$ spectral subspace of $\acxcg$ is the image of $P_n$, so we will show that Im$\ P_i=B_i$, for $i=0,1$. It is clear that the restriction of $P_i$ to $B_i$ is the identity map, so Im$\ P_i\supset B_i$, for $i=0,1$.

Let $\phi\in C_c(G,\acx)$ be such that for some $n > 1$ 
\[\phi(t)=i_X(x_1(t))i_X(x_2(t))\cdots i_X(x_n(t)),\]
where $x_i(t)\in X$ for all $t\in\supp \ \phi$.
Then $\gamma_{\lambda}(\phi)=\lambda^n\phi$, and $P_i(\phi)=0$, for $i=0,1$.
Analogously, $P_i(\phi)=0$ for $i=0,1$  if $\phi\in C_c(G,\acx)$ is such that  for some  $n>0$
\[\phi(t)=i_X(x_1(t))^*i_X(x_2(t))^*\cdots i_X(x_n(t))^* \] for all $t\in\supp \ \phi$.

Let ${\mathcal C}$ denote the dense *-subalgebra of $\acx$ generated by $\{i_A(A),\ i_X(X)\}$. Given $\phi\in C_c(G,\acx)$ and $\epsilon >0$, let $U$ be a precompact open set containing $\supp\ \phi$, and let $\epsilon'=\epsilon/\mu(U)$, $\mu$ being Haar measure. 

For each $t\in \supp\ \phi$ choose $c_t\in {\mathcal C}$ and a neighborhood $N_t\subset U$ of $t$ such that $\|\phi(s)-c_t\|<\epsilon '$ for all $s\in N_t$.
 Let $\{N_{t_i}\}$ be a finite subcovering of $\supp\ \phi$ and $\{h_i\}$ a partition of unity subordinate to it. Then $\|\phi-\sum_ih_ic_{t_i}\|_{{\scriptscriptstyle{\acxcg}}}\leq \|\phi-\sum_ih_ic_{t_i}\|_{{\scriptscriptstyle{L^1(G,\acx)}}}< \epsilon$.

Therefore $P_i(\phi)\in B_i$, and  
\[P_i(\acxcg)\subset \overline P_i(C_c(G,\acx))\subset B_i,\]
for $i=0,1$.

We next show that $\gamma$ is semi-saturated, i.e that $\acxcg$ is   generated as a \cstar-algebra by $B_0$ and $B_1$. As above, any function in $C_c(G,\acx)$ can be approximated by finite sums $\sum f_i$, where  either $f_i$ or $f_i^*$ is of the form 
\[f(t)=h(t)i_A(x_0)i_X(x_1)i_X(x_2)\dots i_X(x_n),\] 
for $h\in C_c(G)$, $n\geq 0$, $x_0\in A$, and $x_1,x_2\dots x_n\in X$. Therefore it suffices to show that these maps  belong to $C^*(B_0,B_1)$. We show this by induction on $n$.  Note that the result holds for $n=0,1$. 

Given   $\epsilon >0$ and $f$ as above for $n>1$, let $V$ be a neighborhood of $e$ in $G$ such that $\|h(t)x_n-h(s^{-1}t)(\alpha_X)_s(x_n)\|<\epsilon$ for all $t\in \supp \ h$ and $s\in V$. Let $\lambda\in C_c(G)$ be a positive function such that $\supp \ \lambda\subset V$ and $\int_G\lambda=1$.
We now set $k(t)=h(t)i_X(x_n)$ and $g(t)=\lambda(t)y$, where $y=i_A(x_0)i_X(x_1)i_X(x_2)\dots i_X(x_{n-1})$, so that $g,k\in \cstar (B_0,B_1)$.

Then
\[\|f(t)-(gk)(t)\|=\|\int_G\lambda(s)y[h(t)i_X(x_n)-h(s^{-1}t)i_X((\alpha_X)_s(x_n)]ds\|< \|y\|\epsilon.\]
Thus $C_c(G,\acx)\subset \cstar(B_0,B_1)$, and $\gamma$ is semi-saturated. 

We now show that $B_0$ is isomorphic to $A\rtimes_{\alpha_{{\scriptscriptstyle A}}}G$. Let $j_{\acx}$ and $j_G$ denote the canonical inclusions of $\acx$ and $G$ in the multiplier algebra of $\acxcg$, respectively. The non-degenerate $^*$-homomorphism $j_{\acx}\circ i_A$ is non-degenerate and the pair
\[(j_{\acx}\circ i_A,j_G):(A,G)\longrightarrow M(\acxcg)\] is covariant for the system $(A,G,\alpha_{{\scriptscriptstyle A}})$. Therefore it  induces, as in \cite[Prop. 2]{raeb}, a $^*$-homomorphism  $J:A\rtimes_{\alpha_{{\scriptscriptstyle A}}} G\longrightarrow \acxcg$ such that $(J\phi)(t)=i_A(\phi(t))$ for all $\phi\in C_c(G,A)$. This shows that the image of $J$ is $B_0$, and it remains to prove that $J$ is one-to-one.

     Let $\pi$ be a   non-degenerate faithful representation of $A$ on a Hilbert space $H$, and let $ {\tilde {\pi}}:=\Lambda\otimes \id_H$ be the representation  of  $\acx$ on $ E_X\otimes_{\pi}H$ induced by $\pi$ as in Definition \ref{indrep}. Denote by  $\theta$  the representation of $\acxcg$ on $L^2(G, E_X\otimes_{\pi}H)$ induced by  ${\tilde {\pi}}$.

Let $V$ be  the unitary operator  from $A\otimes_\pi H$ to $H$ defined in Remark \ref{rkindrep}. Note that $A\otimes_\pi H\subset  E_X\otimes_\pi H$ is invariant under $\tilde{\pi}(i_A(a))$ for all $a\in A$ and that
\[V\tilde{\pi}(i_A(a))V^*=\pi(a)\mbox{ for all } a\in A,\]
because for $a,b\in A$ and $h\in H$
\[V\tilde{\pi}(i_A(a))V^*(\pi(b)h)=V(ab\otimes h)=\pi(a)(\pi(b)h).\]

Fix  now $\xi\in L^2(G, H)$  and  $\phi\in C_c(G,A)\subset A\rtimes_{\alpha_A} G$. Then $V^*\circ \xi\in L^2(G,A\otimes_\pi H )\subset L^2(G,E_X\otimes_\pi H)$, and  
\[\begin{array}{ll}
\left(\theta_{J(\phi)}(V^*\circ\xi)\right)(r)&=\int_G {\tilde \pi}[\alpha_{r^{-1}}(J(\phi)(t))]V^*[\xi(t^{-1}r)]dt\\
\\ 
& =V^*\big(\int_G \pi[ (\alpha_A)_{r^{-1}}(\phi(t))](\xi(t^{-1}r))\big) dt\\
& \\
&= [V^*\pi_0(\phi)(\xi)](r),
\end{array}\]
for all $r\in G$, where $\pi_0$ denotes the representation of $A\rtimes_{\alpha_A}G$ on $L^2(G, H)$ induced by $\pi$.
Thus,
\[\|J(\phi)\|=\|\theta_{J(\phi)}\|\geq \|\pi_0(\phi)\|=\|\phi\|,\]
which shows that $J$ is isometric.

Now, $B_1$ carries  a natural structure of Hilbert \cstar-bimodule over $B_0$,  both the left and the right action consisting of  multiplication, and the  $B_0$-valued inner products being given by $\langle b,c\rangle_L=bc^*$, $\langle b,c\rangle_R=b^*c$. On the other hand, the map $f\mapsto i_X\circ f$ and the isomorphism $J$ above identify $C_c(G,X)$ with  $C_c(G,i_X(X))\subset B_1$ and $B_0$ with $A\rtimes_{\alpha_{\scriptscriptstyle A}} G$, respectively.
It only remains to show that under these identifications $(B_0,B_1)$ and $ (A\rtimes_{\alpha_{\scriptscriptstyle A}} G, Y)$  agree as Hilbert \cstar-bimodules.  For  $f,g\in C_c(G,X)$ we have
\[\begin{array}{ll}
i_A[\langle f,g \rangle_L(t)]&=(i_X\circ f)(i_X\circ g)^*(t)\\
\\
&=\int_G (i_X\circ f)(u)\alpha_u[(i_X\circ g)^*(u^{-1}t)]du\\
\\
&=\int_G i_X(f(u))[i_X(\Delta(t^{-1}u)(\alpha_X)_t(g(t^{-1}u)))]^* du\\
\\
&=i_A\big(\int_G  \Delta(t^{-1}u)\langle f(u),(\alpha_X)_t[g(t^{-1}u)]\rangle_L du \big).
\end{array}\]
The remaining cases are shown by means of similar  computations.
\end{proof}
\section{The duality theorem}

In this section we make use of the results in Proposition \ref{dcp} in order to discuss the crossed product by the dual action on the crossed product by a Hilbert \cstar-bimodule. We establish in Theorem \ref{takai} conditions that ensure the Morita equivalence between $A$ and the double crossed product $\acx\rtimes_\delta S^1$.  
\begin{prop}\label{czero}
Let $X$ be a Hilbert \cstar-bimodule over $A$, and let $C_0(\ZZ,X)$ be the Hilbert \cstar-bimodule over $C_0(\ZZ, A)$ defined by
\begin{align}
(\phi f)(n)=\phi(n)f(n),\quad (f\phi)(n)=f(n)\phi(n-1),\label{E:bimod1}\\
\langle f,g\lran (n)=\langle f(n),g(n)\lran,\quad \langle f,g\rran (n)=\langle f(n+1),g(n+1)\rran,\label{E:bimod2}
\end{align}
for $\phi\in C_0(\ZZ, A)$ and $f,g\in C_0(\ZZ,X)$.

Let $\delta$ denote the dual action of $S^1$ on  $\acx$.
Then there is an isomorphism 
\[I:(\acx)\rtimes_{\delta} S^1\longrightarrow C_0(\ZZ, A)\rtimes C_0(\ZZ,X)\]
such that
\[(I(\phi))(n)=\int_{S^1} \lambda^ni_A^{-1}(\phi(\lambda))d\lambda \mbox{ and } (I(f))(n)=\int_{S^1} \lambda^{n-1}i_X^{-1}(f(\lambda))d\lambda,\]
for $\phi\in C(S^1,i_A(A))$ and $f\in C(S^1,i_X(X))$.

\end{prop}
\begin{proof}
In the notation of Remark \ref{covact},  $\delta$ is the action induced by $(\delta_{{\scriptscriptstyle A}},\delta_{{\scriptscriptstyle X}})$, where, for  $\lambda\in S^1$,  $(\delta_{{\scriptscriptstyle A}})_{\lambda}$ is the identity and $(\delta_{{\scriptscriptstyle X}})_{\lambda}$    is multiplication by $\lambda$.

It follows from Proposition \ref{dcp} that
\[(\acx)\rtimes_{\delta} S^1\simeq (A\rtimes_{\id}S^1)\rtimes Y,\]
$Y$ being the completion of $C(S^1,X)$ with the norm coming from the  $A\rtimes _{\id} S^1$-Hilbert \cstar-bimodule structure given by:
\begin{gather}
(\phi f)(\mu)=\int_{S^1}\lambda\phi(\lambda)f(\lambda^{-1}\mu)d\lambda,\quad (f\phi)(\mu)=\int_{S^1}f(\lambda)\phi(\lambda^{-1}\mu)d\lambda, \label{E:dcp1}\\
\langle f,g\lran (\mu)=\int_{S^1}\langle f(\lambda),\mu g(\mu^{-1}\lambda)\lran d\lambda,\quad \langle f,g\rran (\mu)=\int_{S^1}\langle f(\lambda), g(\lambda\mu)\rran d\lambda.\label{E:dcp2}
\end{gather}
 Let $J_{\sss A}$ denote the isomorphism   $J_{{\sss A}}:A\rtimes_{\id}S^1\longrightarrow C_0(\ZZ,A)$ given by  
\[(J_{{\sss A}}\phi)(n)=\int_{S^1} \lambda^n\phi(\lambda)d\lambda, \]
for $\phi\in C(S^1,A)$, and define $J_{{\sss Y}}$ on $C(S^1,X)$ by 
\[ \ (J_{{\sss Y}}f)(n)=\int_{S^1} \lambda^{n-1}f(\lambda)d\lambda.\label{E:fourier}\]
Note that for $f,g\in C(S^1,X)$ we have
\[\begin{array}{ll}  
J_{\sss A}(\langle f,g \rran )(n)&= \int_{S^1\times S^1} \lambda^n\langle f(\mu),g(\mu\lambda)\rran d\mu\  d\lambda\\
&\\
&=\int_{S^1\times S^1} \mu^{-n}\lambda^n\langle f(\mu),g(\lambda)\rran d\mu\  d\lambda\\
&\\
&=\langle (J_{\sss Y}(f))(n+1), (J_{\sss Y}(g))(n+1)\rran .
\end{array}\]
Therefore
\[\lim_n\|J_{{\sss Y}}(f)(n)\|^2=\lim_n\|J_{{\sss A}}(\langle f,f\rran)(n)\|=0\]
and 
\[\|f\|^2=\|\langle f,f \rran\|=\|J_{\sss A}(\langle f,f\rran)\|=\|\langle J_{\sss Y}(f),J_Y(f)\rran \|=\|J_{\sss Y}(f)\|^2,\]
which shows that $ J_{{\sss Y}}$ extends to $ J_{{\sss Y}}:Y\rightarrow C_0(\ZZ,X)$. In fact, $J_{{\sss Y}}(Y)$ is all of $C_0(\ZZ,X)$, since the Hilbert \cstar-bimodule norm in $C_0(\ZZ,X)$ is the supremum norm, $J_{\sss Y}$ is linear and isometric, and $J_{\sss Y}(f_{x,k})=x\delta_{1-k}$ for $x\in X$,  $k\in \ZZ$,   $f_{x,k}(\lambda)=\lambda^kx$.

Routine computations similar to  those above  show that ($J_{{\sss A}}$, $J_{{\sss Y}}$) is an isomorphism of   Hilbert \cstar-bimodules between $(A\rtimes_{\id} S^1,Y)$, as  defined by Equations (\ref{E:dcp1}) and  (\ref{E:dcp2}), and ($C_0(\ZZ, A), C_0(\ZZ,X))$), as defined by  equations (\ref{E:bimod1}) and (\ref{E:bimod2}).
Finally, it follows from Remark \ref{covmaps} that $J_{\sss A}\rtimes J_{\sss X}$ is an isomorphism from $(A\rtimes_{\id}S^1)\rtimes Y$ to  $ C_0(\ZZ, A)\rtimes C_0(\ZZ,X)$; the isomorphism $I$ is obtained by composing $J_A\rtimes J_X$ with the isomorphism in Proposition \ref{dcp}.

\end{proof}

\bprop
\label{pi}
Let $X$ be a Hilbert \cstar-bimodule over $A$, and let 
\[(\pi_0,\pi_1): (C_0(\ZZ,A),C_0(\ZZ,X))\longrightarrow \cl(E_X)\]
be defined by
\[[(\pi_0\phi)(\eta)](n)=\phi(n)\eta(n)\quad [(\pi_1f)(\eta)](n)=T^{n-1}_{f(n)}((\eta(n-1)),\]
for $\phi\in C_0(\ZZ,A)$ and $f\in C_0(\ZZ,X)$, where $T^{n-1}_{f(n)}$ and $E_X$ are as in Definition \ref{TnEx}. 
Then
\begin{enumerate}

\item The pair $(\pi_0,\pi_1)$ is a representation, so  it induces a $*$-homomorphism
$\pi:C_0(\ZZ,A)\rtimes C_0(\ZZ,X)\longrightarrow \cl(E_X),\ and \ \im\ \pi$ is the \cstar-subalgebra of $\cl(E_X)$  generated by $\{T^n_x,L^n_a:n\in \ZZ, x\in X, a\in A\}$, where $L^n_a$ is as in Definition \ref{TnEx}, and the set above is viewed as a subset of $\cl(E_X)$ as in Remark \ref{included}.

\item The map $\pi$ is injective if and only  if so are the homomorphisms $A\mapsto \cl(X_{{\sss A}})$ and $A\mapsto \cl(_{{\sss A}}X)$ induced by the left and the right action of $A$ on $X$, respectively.

\item The image of $\pi$ contains ${\mathcal K}(E_X)$, and it is ${\mathcal K}(E_X)$ when $X$ is  full both as a left and a right Hilbert \cstar-module over $A$.

\end{enumerate}
\eprop

\begin{proof}

\end{proof} (1)  In what follows  $\phi\in C_0(\ZZ,A)$, $f,g\in C_0(\ZZ,X)$, $\xi,\eta\in E_X$, and $n\in \ZZ$. 
By virtue of equations (\ref{E:tn1}), (\ref{E:tn2}),  (\ref{E:bimod1}), and (\ref{E:bimod2}), we have
\[\begin{array}{ll}
[\pi_1(f\phi)(\eta)](n)&=T^{n-1}_{f(n)\phi(n-1)}(\eta(n-1))\\
&\\
&=T^{n-1}_{f(n)}(\phi(n-1)\eta(n-1))\\
&\\
&=[(\pi_1f)(\pi_0\phi)(\eta)](n)].
\end{array}\]
Note that, since
\[\begin{array}{ll}
\langle(\pi_1f)(\eta),\xi\rangle &= \sum_n\langle T^{n-1}_{f(n)}(\eta(n-1)),\xi(n)\rran\\
&\\ 
&= \sum_n\langle \eta(n),(T^{n}_{f(n+1)})^*(\xi(n+1))\rran,
\end{array}\]
we have
\begin{align}
[(\pi_1f)^*(\xi)](n)=\big(T^n_{f(n+1)}\big)^*(\xi(n+1)).\label{E:pifstar}
\end{align}
Therefore
\[\begin{array}{ll}[\pi_0\big(\langle f,g \rran\big)(\eta)](n)&=\langle f,g \rran(n)\eta(n)\\
&\\
&=\langle f(n+1),g(n+1)\rran \eta(n)\\
&\\
&=[\big(T^n_{f(n+1)}\big)^* T^n_{g(n+1)}](\eta(n))\\
&\\
&=[((\pi_1f)^*(\pi_1g))(\eta)](n).
\end{array}\]
The remaining properties are checked in a similar fashion. The last statement follows from the fact that
\begin{align}
\pi_0(a\delta_n)=L^n_a,\ \pi_1(x\delta_{n+1})=T^n_x,
\end{align}
for $a\in A,\ x\in X, n\in \ZZ$.

(2) Let $U_\lambda$ be the unitary operator defined in Proposition \ref{action}. Conjugation by $U_\lambda$ yields, as in Proposition \ref{action}, a strongly continuous action $\beta$ of $S^1$ on the image of $\pi$, and $\pi $ is covariant for the dual action and $\beta$.

Therefore, by \cite[2.9]{cact}, $\pi$ is injective if and only if so is $\pi_0$. Clearly the left and right actions of $A$ on $X$ must be faithful when $\pi$ is injective because they correspond to $\pi_0$ on maps supported in $\{1\}$ and $\{-1\}$, respectively.

On the other hand, if the left and right actions of $A$ on $X$ are injective, then so are the left and right actions of $A$ on $X^{\otimes n}$ for all $n\in \ZZ$. This is shown by induction on $n$ for positive values of $n$,   since $a (x_1\otimes\dots\otimes x_n)=0$ for all $x_i\in X$ implies that
\begin{align*}
0=\langle|ax_1|(x_2\otimes\dots\otimes x_n), |ax_1|(x_2\otimes\dots\otimes x_n)\rangle_R,  
\end{align*}
for all $x_i\in X$, where  $|ax_1|=(\langle ax_1,ax_1\rran)^{1/2}$. The equality above implies that $|ax_1|=0$ for all $x_1\in X$ and, consequently, that  $a=0$.  A similar reasoning with the left inner product yields the proof for the right action of $A$. It is apparent that the statement always holds for $n=0$. Finally,  the case $n<0$ is taken care of by applying the results above  to  the dual bimodule $\tilde{X}$.

Assume now that the left and right actions of $A$ on $X$ are faithful, and let  $\phi\in \ker \ \pi_0$. Then 
\[0=((\pi_0\phi)\eta) =\phi(n)\eta(n).\]
for all  $\eta\in E_X $ and $n\in \ZZ$. It  follows from the remarks above that $\phi=0$.

(3) First note that $\kex$ is generated as a \cstar-algebra by the set 
\[\{\theta_{\eta\delta_k,\xi\delta_l}:\ k,l\in \ZZ,\ k\geq l,\ \eta=x_1\otimes \dots \otimes x_k, \xi=y_1\otimes \dots \otimes y_l\}.\]

Notice also that for a positive integer $n$, $x\in X$, $u\in X^{\otimes n-1}$, and $\eta\in E_X$ we have
\[\theta_{(x\otimes u)\delta_n, \eta} =T^{n-1}_{x}\theta_{ u\delta_{n-1},\eta}=\pi_1(x\delta_n)\theta_{ u\delta_{n-1},\eta}\]
Similarly, for a positive integer $m$, $\tilde{y}\in \tilde{X}$, $\eta\in E_X$, and $v\in {\tilde{X}}^{\otimes m-1}$,
\[\theta_{\eta,{(\tilde{y}}\otimes v)\delta_{-m}}=\theta_{\eta, v\delta_{-m+1}}\pi_1(\tilde{y}\delta_{-m+1}) ,\]

Finally, since  $\theta_{a,b}=\pi_0(ab^*\delta_0)$ for all $a,b\in A$, and 
$\theta_{\eta,\xi}=\theta^*_{\xi,\eta}$ for all $\eta,\xi\in E_X$, the inclusion ${\mathcal{K}}(E_X)\subset \im\  \pi$ follows.

When $X$ is full on the left and the right,  so are the bimodules $X^{\otimes n}$ for all $n\in \ZZ$, and consequently  $A$ acts  by compact operators on $X^{\otimes n}$ for all  $n\in \ZZ$.
Besides, $T^n_x$ is compact for all $n\in \ZZ$ and $x\in X$: for $n\geq 0$  approximate $x$ by \[x'=x_0\sum_{i=1}^N \langle u_i,v_i \rangle^A_L,\]
for appropriate $u_i,v_i\in X^{\otimes n}$. Then $T^n_x$ gets approximated by 
\[T^n_{x'}=\sum_{i=1}^N\theta_{x_0\otimes u_i,v_i }\in {\mathcal K}(X^{\otimes n},X^{\otimes n+1} ) .\]
Therefore $(T^n_x)^*$ is compact when $n\geq 0$, which shows, by  equation \ref{E:tnegn}, that so is   $T^n_x$ for negative values of $n$. It follows from  Remark \ref{included} that  the images  of both $\pi_0$ and $\pi_1$, and consequently that of $\pi$,  are contained in ${\mathcal K}(E_X)$.

\begin{rk}\label{bidual}
Let $A, X,$ and  $\pi$ be as in Proposition \ref{pi} and identify $(\acx)\rtimes_{\delta} S^1$ and $ C_0(\ZZ, A)\rtimes C_0(\ZZ,X)$ through the isomorphism $I$ in Proposition \ref{czero}. Then, when  $\pi$ is injective,  
the bidual action  of $\ZZ$ on $\im\ \pi$ becomes the automorphism $\sigma$  given by
\[\sigma(L^n_a)=L^{n-1}_a,\ \sigma(T^n_x)=T^{n-1}_x,\]
for all $a\in A$, $x\in X$, and $n\in \ZZ$. (Notice that by part (1) of Proposition \ref{pi} $\sigma$ is determined by the equations above)
\end{rk}
\begin{proof}
In the notation of Propositions \ref{czero} and \ref{pi} we have
\[L^n_a=\pi(a\delta_n)=\pi \circ I(f_{n,a}),\ \mbox{ where }\ f_{n,a}(\lambda)=\lambda^{-n}a,\]
for all $a\in A$ and $n\in \ZZ$.

Then
\[\sigma(L^n_a)=\pi\circ I(\hat{\delta}(f_{n,a}))=\pi \circ I(f_{n-1,a})= L^{n-1}_a.\]
Analogously, $T^n_x=\pi\circ I (g_{n,x})\ \mbox{,  for }\  g_{n,x}(\lambda)=\lambda^{-n}x,$ $x\in X$,  and $n\in \ZZ$,
 and  $ \sigma(T^n_x)=T^{n-1}_x.$
\end{proof}

\begin{thm}
\label{takai}
Let $X$ be an $A-A$ Morita equivalence bimodule,  that is,   a Hilbert \cstar-bimodule over $A$ that is full both on the left and the right. Let $\delta$ denote the dual action of $S^1$ on $\acx$. Then the crossed product $(\acx)\rtimes_\delta S^1$ is Morita equivalent to $A$.

Namely, there is an isomorphism between  $(\acx)\rtimes_\delta S^1$ and  ${\mathcal K}(E_X)$ through which the bidual action $\hat{\delta}$ becomes the action induced by the automorphism  $\sigma$  given   by
\[\sigma(L^n_a)=L^{n-1}_a,\ \sigma(T^n_a)=T^{n-1}_x,\]
for all $a\in A$, $x\in X$, and $n\in \ZZ$.

\end{thm}
\begin{proof}
It suffices to show the second statement, since $E_X$ is a Morita equivalence bimodule between $\ck(E_X)$ and $A$.
In view of Proposition \ref{czero}, Proposition \ref{pi},  and Remark \ref{bidual}, it only remains to notice that the left and right actions of $A$ on $X$ are faithful because $X$ is full.
\end{proof}

\begin{clly}
If $X$ is  a Morita equivalence bimodule over $A$, then $\acx$ is Morita equivalent to the crossed product (by an automorphism) $\ck(E_X)\rtimes_{\sigma}\ZZ$, where
$\sigma$ is the automorphism defined in Remark \ref{bidual}.
\end{clly}
\begin{proof}
By Takai duality (\cite{takai}), $\acx$ is Morita equivalent to $\acx\rtimes_{\delta}S^1\rtimes_{\hat{\delta}}\ZZ$, $\delta$ being the dual action on $\acx$. Now, by Theorem \ref{takai}, $\acx\rtimes_{\delta}S^1\rtimes_{\hat{\delta}}\ZZ$
  is isomorphic to 
$\ck(E_X)\rtimes_{\sigma}\ZZ$.
\end{proof}

\begin{rk}
\label{auto}
\rm{Let $\alpha\in \aut(A)$, so that $A\rtimes_\alpha\ZZ$ is isomorphic to $A\rtimes X$  for  $X=A_{\alpha}$ as in Remark \ref{aalpha}. By identifying $A^{\otimes n}$ and $A$ through the isomorphism  $I_n$  in Remark \ref{aalpha}, we get an isomorphism $I$ of right Hilbert modules from $E_X$ to $l^2(\ZZ)\otimes A$. Namely, $(I\eta)(n)=I_n(\eta(n))$.

This isomorphism yields, in turn, an identification between ${\mathcal K}(E_X)$ and ${\mathcal K}\otimes A$ which,  by virtue of equations (\ref{E:in1}) and (\ref{E:in2}), maps $L^n_a$ and $T_x^n$ to $E_{nn}\otimes \alpha^{-n}(a)$ and $E_{n+1,n}\otimes \alpha^{-(n+1)}(x)$ , respectively, where, as usual, $E_{ij}$ denotes the matrix having a 1 at its $ij$-entry and all other entries 0.
In this setting the automorphism $\sigma$ in Theorem \ref{takai} becomes, as in \cite{takai}, Ad$\rho\otimes \alpha$, $\rho$ being translation by 1.}
\end{rk}

 \end{document}